\documentclass[12pt,leqno]{article}
\usepackage{a4wide,epsf,latexsym}
\usepackage[centertags]{amsmath}
\usepackage{graphicx}
\numberwithin{equation}{section}

\newcommand{\ve}{\varepsilon}

\begin{document}

\title{Remarks on boundary layers in singularly perturbed Caputo fractional boundary value
problems}

\author{Hans-G. Roos}

\date{}

\maketitle

\begin{abstract}
Almost nothing is known about the layer structure of solutions to
singularly perturbed Caputo fractional boundary value problems. We discuss
simple convection-diffusion and reaction-diffusion problems.
\end{abstract}


\section{Introduction}
Fractional differential problems became more and more popular in the last years. For an introduction into fractional differential equations
see \cite{Di10}, for instance.

But in the singularly perturbed case we know not much
about the existence of layers. In the case of Caputo fractional boundary value problems we discuss two simple examples which
show that fractional problems are different from standard boundary value problems.

\section{Convection-diffusion}
Consider the singularly perturbed boundary value problem
	\begin{equation}
		-\ve D^{2-\alpha}_*u-u' =f,\qquad
		u(0)=u(1)=0
	\end{equation}
with $0<\ve<<1$, $\alpha\in (0,1)$ and the Caputo derivative $D^{2-\alpha}_*$. For simplicity
we take $f=-1$ and often choose  $\alpha=1/2$ to study
\begin{equation}
		-\ve D^{3/2}_*u-u' =-1,\qquad
		u(0)=u(1)=0.
	\end{equation}
In the classical case  $\alpha=0$ we can decompose $u$ into $u=U+V$, where $U$ solves the reduced
problem
\begin{equation}
		-U' =-1,\qquad
		U(1)=0,
	\end{equation}
while the layer component $V(\xi)$ with $\xi=x/\ve$ solves
\begin{equation}
		D^{2}V+DV =0,\qquad
		V(0)=-U(0),\quad V(1/\ve)=0.
	\end{equation}
Thus we get the exponentially decaying layer component
\[
V=-U(0)\frac{e^{-x/\ve}-e^{-1/\ve}}{1-e^{-1/\ve}}\approx -U(0)e^{-x/\ve}.
\]
{\it Surprisingly, the case  $\alpha\in (0,1)$ is very different from the case $\alpha=0$.}

Let us in this case again decompose
\[
u=U+V_\alpha.
\]
With $\xi=\frac{x}{\ve^{1/(1-\alpha)}}$ we get
\begin{equation}\label{bvp}
		D^{2-\alpha}_*V\alpha+DV_\alpha =0,\qquad
		V_\alpha(0)=-U(0),\quad V_\alpha(1/\ve^{1/(1-\alpha)})=0.
	\end{equation}
In our special case we have $U(0)=-1$.

Next we solve instead the initial value problem
\begin{equation}\label{ivp}
		D^{2-\alpha}_*V\alpha+DV_\alpha =0,\qquad
		V_\alpha(0)=1,\quad DV_\alpha(0)=\theta
	\end{equation}
and try to determine $\theta$ in such a way that the solution of \eqref{ivp} solves the
boundary value problem \eqref{bvp} as well. One possibility for solving  \eqref{ivp} consists
in the application of the Laplace transform. We get
\[
\hat V_\alpha(s)=\frac{s^\alpha+s+\theta}{s^{1+\alpha}(1+s^{1-\alpha})}.
\]
Now it is possible \cite{Be74} from the asymptotic behavior of the function $\hat V_\alpha$ to conclude
the asymptotic behavior of the original function. We observe:\\
{\it $V_\alpha$ is exponentially decaying if and only if $\alpha=0$ and $\theta=-1$}.\\
Otherwise we have $V_\alpha\sim \xi^\alpha$ if $\xi\mapsto \infty$.

Alternatively we solve \eqref{ivp} for $\alpha=1/2$ by power series. Setting
\[
V_{1/2}=1+\theta\,V_{1/2}^*
\]
we obtain
\[
V_{1/2}^*=\sum_{k=1}^\infty\frac{\xi^k}{k!}-\sum_{k=1}^\infty\frac{\xi^{k+1/2}}{\Gamma(k+3/2)}.
\]
Using Mittag-Leffler functions \cite{HMS09} one can write
\begin{equation}\label{M-L}
V_{1/2}(\xi)=1+\theta\,\xi\,E_{1/2,2}(-\xi^{1/2})
\end{equation}
(this also follows from the Laplace-transform of $V_{1/2}$).

Consequently, we obtain for the layer correction with $\theta=-1/V_\alpha^*(1/\ve^2)$ in the
case $\alpha=1/2$
\begin{equation}
V_{1/2}(x)=1-\frac{V^*_{1/2}(x/\ve^2)}{V^*_{1/2}(1/\ve^2)}.
\end{equation}
From the asymptotic behavior of the Mittag-Leffler functions we conclude
\[
V^*_{1/2}(1/\ve^2)=O(1/\ve)\quad{\rm and}\,\,\, \theta=O(\ve).
\]
For fixed $x_0$ we obtain
\[
\lim_{\ve\to 0}V_{1/2}(x_0)=1-\lim_{\ve \to 0}\frac{V^*_{1/2}(x_0/\ve^2)}{V^*_{1/2}(1/\ve^2)}=
    1-\frac{\sqrt{x_0}/\ve}{1/\ve}=1-\sqrt{x_0}.
\]
Thus, $V_{1/2}$ has nothing to do with a classical boundary layer function.
See Figure 2: left the solution, right the layer correction.

\section{Reaction-diffusion}
Consider the boundary value problem
\begin{equation}
		-\ve D^{2-\alpha}_*u+u =f,\qquad
		u(0)=u(1)=0
	\end{equation}
with $0<\ve<<1$, $\alpha\in (0,1)$ and the Caputo derivative $D^{2-\alpha}_*$. For simplicity
we choose $f=-1$ and mainly discuss the case $\alpha=1/2$.

We decompose the solution into $u=-1+V_\alpha^0++V_\alpha^1$. Then $V_\alpha^0$ and $V_\alpha^1$ satisfy
the homogeneous equation, moreover
\[
 V_\alpha^0(0)=1,\,\,V_\alpha^0(1)=\mu;\quad V_\alpha^1(0)=0,\,\,V_\alpha^1(1)=1-\mu.
\]
Later we will fix $\mu$ and show that it is small.

{\it Are $V_\alpha^0$ and $V_\alpha^1$ boundary layer functions?}

At $x=0$ we introduce the local variable $\xi=x/\ve^{1/(2-\alpha)}$. Thus, for instance for $\alpha=1/2$ we obtain
for $V_\alpha^0$
\[
  -\ve D^{3/2}_*V+V =0,\qquad
		V(0)=1,\,\,V(1/\ve^{2/3})=0.
\]
Again we replace the boundary value problem by an initial value problem with $DV(0)=\theta.$ Its
Laplace transform yields
\[
\hat V^0_\alpha(s)=\frac{s+\theta}{s^{\alpha}(s^{2-\alpha}-1)}.
\]
Consequently, for $\theta =-1$ there exists a solution $ V^0_\alpha(\xi)$  that decays as $\xi^{\alpha-1}$.
We set $\mu=V^0_\alpha(\xi)|_{x=1}$ and get $\mu=O(\ve^{(1-\alpha)/(2-\alpha)})$.

For $V_\alpha^1$ we solve the correspondent initial value problem with $DV(0)=\theta$ using the
Laplace transform. We get
\[
\hat V_\alpha^1(\theta)=\frac{\ve\theta}{s^{\alpha}(\ve s^{2-\alpha}-1)}.
\]
For instance for $\alpha=1/2$, the nominator has a zero at $s=1/\ve^{2/3}$, thus the original function
is characterized by $ V\sim \exp(x/\ve^{2/3})$. Choosing $\ve\theta=\exp(-1/\ve^{2/3})$, we observe
that $ V_\alpha^1$ behaves like a typical exponentially decaying boundary layer function, i.e.,
\[
     V_\alpha^1(x)\sim (1-\mu)\exp(-(1-x)/\ve^{2/3}).
\]
See Figure 1: left the solution, right the layer correction with different layers at $x=0$ and
$x=1$.

\begin{figure}[htb]
	\centering
  \includegraphics[width=0.7\textwidth,height=0.4\textwidth]{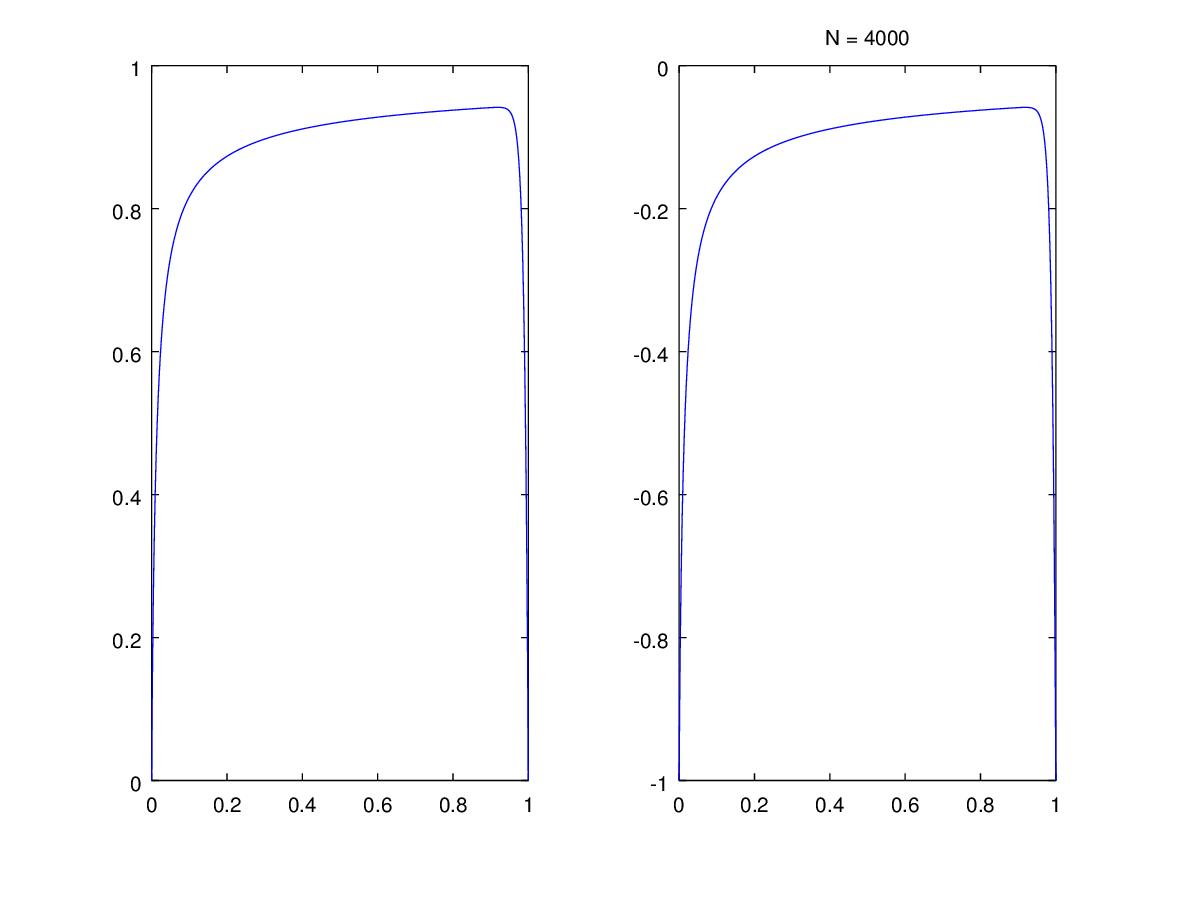}
	\caption{Reaction-Diffusion}
	\label{fig1}
\end{figure}

\begin{figure}[htb]
	\centering
  \includegraphics[width=0.7\textwidth,height=0.4\textwidth]{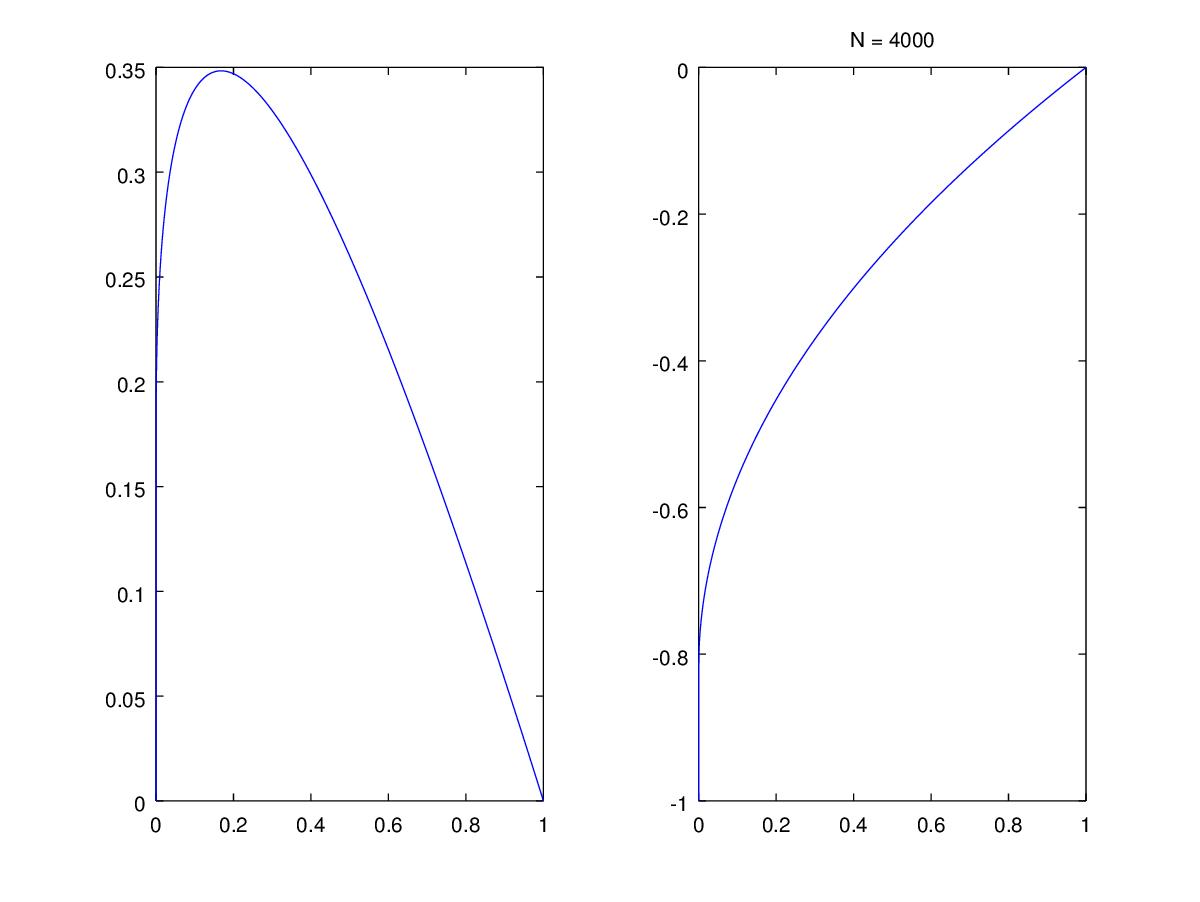}
	\caption{Convection-Diffusion}
	\label{fig2}
\end{figure}


{\bf Acknowledgement} {\it Thanks to Lars Ludwig for generating the two pictures.}

\end{document}